\documentclass[12pt]{article}

\usepackage{amsmath,amsfonts,amssymb,amsbsy}

\usepackage{color}

\newcounter{remark}

\newcounter{sect1}
\newcounter{subsect1}

\newfont{\gotikai}{eufm10}

\newcommand {\beq}{\begin{equation}}
\newcommand {\eeq}{\end{equation}}

\newtheorem{Theorem}{\indent Theorem}
\newtheorem{Lemma}{\indent Lemma}

\def\({\left(}
\def\){\right)}

\def\({\left(}
\def\){\right)}

\def\[{\left[}
\def\]{\right]}

\def\|{\left|}
\def\|{\right|}

\def\0{{\boldsymbol{0}}}

\def\DD{P\ r\ o\ o\ f}

\def\00{{\boldsymbol{0}}}

\usepackage{authblk}
\usepackage[colorlinks]{hyperref}
\usepackage{hyperref}
\hypersetup{
    colorlinks=true,
    linkcolor=blue,
    filecolor=magenta,      
    urlcolor=cyan,
    pdftitle={Overleaf Example},
    pdfpagemode=FullScreen,
    }

\providecommand{\keywords}[1]
{
  \small	
  \textbf{\textit{Keywords---}} #1
}

\begin{document}
\baselineskip = 5.7mm  
	
 	\title{\textbf{Central Limit Theorem on Symmetric Kullback–Leibler (KL) Divergence}}
	
         \author{Helder Rojas \vspace{-2mm}
        \\\small{Escuela Profesional de Ingeniería Estadística, Universidad Nacional de Ingeniería, Perú}
        \\{Department of Mathematics, Imperial College London, United Kingdom}
         \\ email: \href{mailto:me@somewhere.com}{h.rojas-molina23@imperial.ac.uk} \and \vspace{-4mm} Artem Logachov \vspace{-2mm}
         \\\small{Sobolev Institute of Mathematics, Siberian Branch of the Russian Academy of Science, Russia}
        \\{Novosibirsk State Technical University}
        \\{Siberian State University of Geosystems and Technologies}
         \\ email: \href{mailto:me@somewhere.com} {omboldovskaya@mail.ur} }

	\date{} 

	\maketitle
	\vspace{-10mm}
	\begin{abstract}
         In this paper we provide an asymptotic theory for the symmetric version of the Kullback--Leibler (KL) divergence. We define a estimator for this divergence and study its asymptotic properties. In particular, we prove Law of Large Numbers (LLN) and the convergence to the normal law in the Central Limit Theorem (CLT) using this estimator.
         \end{abstract}
         \keywords{Symmetric KL-Divergence estimation,  Central Limit Theorem (CLT).}
\begin{section}{Introduction} 
In mathematical statistics, the Kullback–Leibler (KL) divergence introduced in \cite{kullback}, also called relative entropy in information theory, is a type of statistical measure used to quantify the dissimilarity between two probability measures. This measure of divergence has been widely used in various fields, such as variational inference \cite{blei2017variational, zhang2018advances}, Bayesian inference \cite{tzikas2008variational, jewson2018principles}, metric learning \cite{ji2020kullback, noh2014bias, suarez2021tutorial}, machine learning \cite{claici2020model, poczos2012nonparametric}, computer vision \cite{cui2016feature, deledalle2017estimation}, physics \cite{granero2018kullback}, 
biology \cite{charzynska2015improvement}, information geometry \cite{belavkin2015asymmetric}, among many other application fields. Let $\mathcal{X}=\{a_1, a_2, \dots,a_r\}$ be a finite countable set, where $r\geq 2$. The probability measures on $\mathcal{X}$ are finite dimensional vectors $\mathbf{p}$ in 
\begin{equation*}
    \mathcal{P}(\mathcal{X})=\bigg\{\mathbf{p}=(p_a)_{a\in\mathcal{X}}\in\mathbb{R}^r: p_a\geq 0 \quad\textrm{and}\quad\sum_{a\in\mathcal{X}}p_a=1\bigg\}.
\end{equation*}
If $\mathbf{p}, \mathbf{q}\in \mathcal{P}(\mathcal{X})$, the KL divergence between $\mathbf{p}$ and $\mathbf{q}$ is defined, according to  \cite{kullback}, as following
\begin{equation*}
\mathcal{D}_{KL}(\mathbf{p}||\mathbf{q})=\sum_{a\in\mathcal{X}}p_a\ln\Big(\frac{p_a}{q_a}\Big).
\end{equation*}
The KL divergence is possibly the most famous and used of a large family of divergences called $f$-Divergence. However, KL divergence is asymmetric, i.e., the values of $\mathcal{D}_{KL}(\mathbf{p}||\mathbf{q})$ and $\mathcal{D}_{KL}(\mathbf{q}||\mathbf{p})$ can be significantly different. This lack of symmetry, in some specific contexts, can be a disadvantage when measuring similarities between probability measures, see e.g. \cite{belavkin2015asymmetric, seghouane2006variants, audenaert2013asymmetry, pinski2015kullback, zeng2023fractal, kaminski2022symmetry, johnson2001symmetrizing}. That is why it is often quite useful to work with a symmetrization of the KL divergence which is defined as
\begin{eqnarray}\nonumber
\mathcal{D}_{KL}^{sym}(\mathbf{p}||\mathbf{q})&:=&\mathcal{D}_{KL}(\mathbf{p}||\mathbf{q})+\mathcal{D}_{KL}(\mathbf{q}||\mathbf{p}),\\ \nonumber
    &=&\sum_{a\in\mathcal{X}}(p_a-q_a)\ln\Big(\frac{p_a}{q_a}\Big).
\end{eqnarray}
This symmetrization, also known as Jeffreys divergence, was introduced and used in \cite{jeffreys1998theory}. For some benefits obtained by including symmetrization in KL divergence, see e.g. \cite{chen2018centroid, andriamanalimanana2019symmetric, domke2021easy}. As a consequence of the desirable properties of $\mathcal{D}_{KL}^{sym}(\mathbf{p}||\mathbf{q})$, interest in estimating its value on a sample basis has been growing \cite{nguyen2017supervised, moreno2003kullback, yao2011symmetric, rojas2023statistical, rojas2023order}. Therefore, we consider it important to establish good estimators for the symmetric KL divergence, and in turn, guarantee the efficiency and precision of their estimates through the study of their respective asymptotic properties. The statistical argument presented above is what motivated and guided our interest in obtaining our main results.

Before we state the main results we need a few definitions. Let $(X,Y)$ be a bivariate random vector defined on the probability space $(\Omega, \mathcal{A}, \mathbf{P})$ with $\Omega=\mathcal{X}\times\{0,1\}$. We denote by $\mathbf{E}$ the expectation with respect to $\mathbf{P}$. Let $\mathbf{p}=(p_{j})_{1\leq j \leq r}$ and $\mathbf{q}=(q_{j})_{1\leq j \leq r}$ two probability measures on $\mathcal{X}$ such that, for any $j=1,\dots, r$, we have that
\begin{equation*}
    p_{j}=\mathbf{P}(X=a_{j}|Y=1)\quad \textrm{and} \quad q_{j}=\mathbf{P}(Y=a_{j}|Y=0).
\end{equation*}
Consider the empirical probability measures $\hat{\mathbf{p}}_n=(\hat{p}_{j,n})_{1\leq j \leq r}$ and $\hat{\mathbf{q}}_n=(\hat{q}_{j,n})_{1\leq j \leq r}$, generated by the sequence i.i.d. random variables $(X_1,Y_1)\dots ,(X_n,Y_n)$ from the bivariate distribution,  $(\mathbf{p},\mathbf{q})$, defined from 
\begin{equation}\label{empirical}
 \hat{p}_{j,n}:=\frac{\sum\limits_{i=1}^n\mathbf{I}(X_i=a_j,Y_i=1)}{\sum\limits_{i=1}^n\mathbf{I}(Y_i=1)}, \quad\quad
\hat{q}_{j,n}:=\frac{\sum\limits_{i=1}^n\mathbf{I}(X_i=a_j,Y_i=0)}{\sum\limits_{i=1}^n\mathbf{I}(Y_i=0)}.
\end{equation}
Additionally, consider the following notation
\begin{equation}\label{empirical2}
\hat{p}_n:=\frac{1}{n}\sum\limits_{i=1}^n\mathbf{I}(Y_i=1), \quad\quad
\hat{q}_n:=\frac{1}{n}\sum\limits_{i=1}^n\mathbf{I}(Y_i=0).
\end{equation}
Based on the Equation \eqref{empirical}, we define an estimator for the symmetric KL divergence divergence $\mathcal{D}_{KL}^{sym}(\mathbf{p}||\mathbf{q})$ which has the form
\begin{equation}\label{estimator}
   \mathcal{D}_{KL}^{sym}(\hat{\mathbf{p}}_n||\hat{\mathbf{
   q}}_n):=\sum\limits_{j=1}^r(\hat{p}_{j,n}-\hat{q}_{j,n})\ln\frac{\hat{p}_{j,n}}{\hat{q}_{j,n}}.
\end{equation}
In this paper, we study the convergence properties, in particular LLN and CLT, of the estimator defined in \eqref{estimator}. To this end, consider the following notation
\begin{eqnarray}\nonumber
\eta_n&:=&\mathcal{D}_{KL}^{sym}(\hat{\mathbf{p}}_n||\hat{\mathbf{
   q}}_n)-\mathcal{D}_{KL}^{sym}(\mathbf{p}||\mathbf{q}),\\ \nonumber
    &=&\sum\limits_{j=1}^r\left((\hat{p}_{j,n}-\hat{q}_{j,n})\ln\frac{\hat{p}_{j,n}}{\hat{q}_{j,n}}
-(p_j-q_j)\ln\frac{p_j}{q_j}\right).
\end{eqnarray}
It is true that there are some works on asymptotic theory for divergence measurements, especially for asymmetric KL divergence, see e.g. \cite{bulinski2021statistical, ba2019divergence, yao2024symmetric, bobkov2019renyi}. However, to the best of our knowledge, there are very few theoretical works about this symmetric version of the KL divergence. In this paper we have been careful with mathematical rigor and consequently we have obtained finer and simpler results than those available in the existing literature.
\end{section}

\begin{section}{Main results} 

\begin{Theorem} (\textbf{Law of Large Numbers}) \label{t.1} The following equality holds
$$
\lim\limits_{n\rightarrow\infty}\eta_n =0 \ \ \ \text{a.s.}
$$
\end{Theorem}
\vspace{1mm}
\DD~  Follows directly from strong law of large numbers for sequences
$\hat{p}_n$, $\hat{q}_n$, $\hat{p}_{j,n}$, $\hat{q}_{j,n}$ defined in \eqref{empirical}  and \eqref{empirical2}, $1\leq j \leq r$.$\Box$
\vspace{1mm}

\begin{Theorem} (\textbf{Central Limit Theorem}) \label{t.2} The following convergence holds true
$$
\sqrt{n}\eta_n\stackrel{d}{\rightarrow}\xi\sim\Phi_{0,\sigma^2},
  \ \text{as} \ n\rightarrow\infty,
$$
where
$$
\begin{aligned}
\sigma^2=\mathbf{E}\Biggl(&
\sum\limits_{j=1}^r\Biggl(\left(\frac{1}{p}\sum\limits_{i=1}^n(\mathbf{I}(X_i=a_j,Y_i=1)-pp_j)
-p_j\sum\limits_{i=1}^n(\mathbf{I}(Y_i=1)-p)\right)\left(1+\ln\frac{p_j}{q_j}-\frac{q_j}{p_j}\right)
\\
&+\left(\frac{1}{q}\sum\limits_{i=1}^n(\mathbf{I}(X_i=a_j,Y_i=0)-qq_j)
-q_j\sum\limits_{i=1}^n(\mathbf{I}(Y_i=0)-q)\right)\left(1+\ln\frac{q_j}{p_j}-\frac{p_j}{q_j}\right)\Biggl)\Biggl)^2.
\end{aligned}
$$
\end{Theorem}
\vspace{3mm}
\DD. Let's first perform the transformations we need
\beq\label{27.11.11}
\begin{aligned}
&\sqrt{n}\sum\limits_{j=1}^r\left((\hat{p}_{j,n}-\hat{q}_{j,n})\ln\frac{\hat{p}_{j,n}}{\hat{q}_{j,n}}
-(p_j-q_j)\ln\frac{p_j}{q_j}\right)
\\=&\sqrt{n}\sum\limits_{j=1}^r\left((\hat{p}_{j,n}-\hat{q}_{j,n})\ln\frac{\hat{p}_{j,n}}{\hat{q}_{j,n}}\pm
(p_j-q_j)\ln\frac{\hat{p}_{j,n}}{\hat{q}_{j,n}}
-(p_j-q_j)\ln\frac{p_j}{q_j}\right)
\\=&\sqrt{n}\sum\limits_{j=1}^r((\hat{p}_{j,n}-p_j)-(\hat{q}_{j,n}-q_j))\ln\frac{\hat{p}_{j,n}}{\hat{q}_{j,n}}
+\sqrt{n}\sum\limits_{j=1}^r\left(p_j-q_j\right)\ln\frac{\hat{p}_{j,n}q_j}{p_j\hat{q}_{j,n}}
\\
= &\sqrt{n}\sum\limits_{j=1}^r((\hat{p}_{j,n}-p_j)-(\hat{q}_{j,n}-q_j))\ln\frac{\hat{p}_{j,n}}{\hat{q}_{j,n}}
+\sqrt{n}\sum\limits_{j=1}^r\left(p_j-q_j\right)\ln\left(1+\frac{\hat{p}_{j,n}q_j-p_j\hat{q}_{j,n}}{p_j\hat{q}_{j,n}}\right).
\end{aligned}
\eeq
Let
$$
y_{j,n}:=\ln\left(1+\frac{\hat{p}_{j,n}q_j-p_j\hat{q}_{j,n}}{p_j\hat{q}_{j,n}}\right).
$$
Note that the Taylor formula for the function $\ln(1+ x)$ with a remainder term in Lagrange form can be written in the following way: for $|x|<1$
\beq\label{11.09.11}
\ln(1+ x)=x-\frac{x^2}{2(1+K_xx)^2},
\eeq
where $K_x\in(0,1)$. Applying (\ref{11.09.11}), we obtain for any $1\leq j \leq r$
$$
\left|y_{j,n}
-\frac{\hat{p}_{j,n}q_j-p_j\hat{q}_{j,n}}{p_j\hat{q}_{j,n}}\right|
\leq2\left(\frac{\hat{p}_{j,n}q_j-p_j\hat{q}_{j,n}}{p_j\hat{q}_{j,n}}\right)^2
\mathbf{I}\left(\left|\frac{\hat{p}_{j,n}q_j-p_j\hat{q}_{j,n}}{p_j\hat{q}_{j,n}}\right|\leq\frac{1}{2}\right)
$$
\beq\label{27.11.1}
+\left|y_{j,n}
-\frac{\hat{p}_{j,n}q_j-p_j\hat{q}_{j,n}}{p_j\hat{q}_{j,n}}\right|
\mathbf{I}\left(\left|\frac{\hat{p}_{j,n}q_j-p_j\hat{q}_{j,n}}{p_j\hat{q}_{j,n}}\right|>\frac{1}{2}\right) \ \ \ \text{a.s.}
\eeq
Using (\ref{27.11.1}) and Lemma~\ref{l.1} (see  (\ref{25.11.3}), (\ref{25.11.5})), for any  $\varepsilon>0$, $1\leq j \leq r$ and some constant $C:=C(\varepsilon,p_1,\dots,p_r,q_1,\dots,q_r,p,q)>0$ we obtain
$$
\mathbf{P}\left(\sqrt{n}\left|y_{j,n}
-\frac{\hat{p}_{j,n}q_j-p_j\hat{q}_{j,n}}{p_j\hat{q}_{j,n}}\right|>\varepsilon\right)
\leq\mathbf{P}\left(2\left(\frac{\hat{p}_{j,n}q_j-p_j\hat{q}_{j,n}}{p_j\hat{q}_{j,n}}\right)^2
\mathbf{I}\left(\left|\frac{\hat{p}_{j,n}q_j-p_j\hat{q}_{j,n}}{p_j\hat{q}_{j,n}}\right|\leq\frac{1}{2}\right)>
\frac{\varepsilon}{2\sqrt{n}}\right)
$$
$$
+\mathbf{P}\left(\sqrt{n}\left|y_{j,n}
-\frac{\hat{p}_{j,n}q_j-p_j\hat{q}_{j,n}}{p_j\hat{q}_{j,n}}\right|
\mathbf{I}\left(\left|\frac{\hat{p}_{j,n}q_j-p_j\hat{q}_{j,n}}{p_j\hat{q}_{j,n}}\right|>\frac{1}{2}\right)>\frac{\varepsilon}{2}\right)
$$
$$
\leq\mathbf{P}\left(\left(\frac{\hat{p}_{j,n}q_j-p_j\hat{q}_{j,n}}{p_j\hat{q}_{j,n}}\right)^2
>\frac{\varepsilon}{4\sqrt{n}}\right)+
\mathbf{P}\left(\left|\frac{\hat{p}_{j,n}q_j-p_j\hat{q}_{j,n}}{p_j\hat{q}_{j,n}}\right|>\frac{1}{2}\right)
$$
$$
\leq2\mathbf{P}\left(\left|\frac{\hat{p}_{j,n}q_j-p_j\hat{q}_{j,n}}{p_j\hat{q}_{j,n}}\right|
>\frac{\min(\sqrt{\varepsilon},1)}{2n^{\frac{1}{4}}}\right)
$$
$$
\leq2\mathbf{P}\left(q_j\left|\hat{p}_{j,n}-p_j\right|+p_j\left|\hat{q}_{j,n}-q_j\right|
>\frac{\min(\sqrt{\varepsilon},1)p_jq_j}{4n^{\frac{1}{4}}},\hat{q}_{j,n}\geq\frac{q_j}{2}\right)
+2\mathbf{P}\left(\hat{q}_{j,n}<\frac{q_j}{2}\right)
$$
$$
\leq 2\mathbf{P}\left(\left|\hat{p}_{j,n}-p_j\right|
>\frac{\min(\sqrt{\varepsilon},1)p_j}{8n^{\frac{1}{4}}}\right)
+2\mathbf{P}\left(\left|\hat{q}_{j,n}-q_j\right|
>\frac{\min(\sqrt{\varepsilon},1)q_j}{8n^{\frac{1}{4}}}\right)
$$
\beq\label{26.11.1}
+2\mathbf{P}\left(|q_j-\hat{q}_{j,n}|>\frac{q_j}{2}\right)\leq 4\exp\{-\sqrt{n}C\}+2\exp\{-n C\}.
\eeq
Denote
$$
z_{j,n}:=\frac{\hat{p}_{j,n}q_j-p_j\hat{q}_{j,n}}{p_j\hat{q}_{j,n}}, \ \ \
v_{j,n}:=\frac{\hat{q}_{j,n}-q_j}{q_j}.
$$
For any $1\leq j \leq r$ we have
$$
\left|z_{j,n}-\frac{\hat{p}_{j,n}q_j-p_j\hat{q}_{j,n}}{p_jq_{j}}\right|
=\frac{1}{p_jq_j}
\left|\frac{\hat{p}_{j,n}q_j-p_j\hat{q}_{j,n}}{1+\frac{\hat{q}_{j,n}-q_j}{q_j}}-(\hat{p}_{j,n}q_j-p_j\hat{q}_{j,n})\right|
$$
$$
=\frac{1}{p_jq_j}
\left|(\hat{p}_{j,n}q_j-p_j\hat{q}_{j,n})\left(-v_{j,n}
+\sum\limits_{k=2}^\infty\left(-v_{j,n}\right)^k\right)\right|
\mathbf{I}\left(\left|v_{j,n}\right|\leq\frac{1}{2}\right)
$$
$$
+\left|z_{j,n}-\frac{\hat{p}_{j,n}q_j-p_j\hat{q}_{j,n}}{p_jq_{j}}\right|\mathbf{I}\left(\left|v_{j,n}\right|\geq\frac{1}{2}\right)
$$
$$
\leq\frac{2}{p_j}|\hat{p}_{j,n}-p_j|\cdot|v_{j,n}|\mathbf{I}\left(\left|v_{j,n}\right|\leq\frac{1}{2}\right)
+\frac{2}{q_j}|\hat{q}_{j,n}-q_j|\cdot|v_{j,n}|\mathbf{I}\left(\left|v_{j,n}\right|\leq\frac{1}{2}\right)
$$
\beq\label{27.11.7}
+\left|z_{j,n}-\frac{\hat{p}_{j,n}q_j-p_j\hat{q}_{j,n}}{p_jq_{j}}\right|\mathbf{I}\left(\left|v_{j,n}\right|\geq\frac{1}{2}\right)
\ \ \  \text{a.s.}
\eeq
Applying (\ref{27.11.7}) and Lemma~\ref{l.1} (see (\ref{25.11.3}), (\ref{25.11.5})), for any $\varepsilon>0$, $1\leq j \leq r$ and for some constant
$C:=C(\varepsilon,p_1,\dots,p_r,q_1,\dots,q_r,p,q)>0$ we obtain
$$
\mathbf{P}\left(\sqrt{n}\left|z_{j,n}-\frac{\hat{p}_{j,n}q_j-p_j\hat{q}_{j,n}}{p_jq_{j}}\right|>\varepsilon\right)
\leq\mathbf{P}\left(|\hat{p}_{j,n}-p_j|\cdot|v_{j,n}|\mathbf{I}\left(\left|v_{j,n}\right|\leq\frac{1}{2}\right)>\frac{\varepsilon p_j}{6\sqrt{n}}\right)
$$
$$
+\mathbf{P}\left(|\hat{q}_{j,n}-q_j|\cdot|v_{j,n}|\mathbf{I}\left(\left|v_{j,n}\right|\leq\frac{1}{2}\right)>\frac{\varepsilon q_j}{6\sqrt{n}}\right)
$$
$$
+\mathbf{P}\left(\sqrt{n}\left|z_{j,n}-\frac{\hat{p}_{j,n}q_j-p_j\hat{q}_{j,n}}{p_jq_{j}}\right|
\mathbf{I}\left(\left|v_{j,n}\right|>\frac{1}{2}\right)>\frac{\varepsilon}{3}\right)
\leq\mathbf{P}\left(|\hat{p}_{j,n}-p_j|>\frac{\sqrt{\varepsilon p_j}}{\sqrt{6}n^{\frac{1}{4}}}\right)
$$
$$
+\mathbf{P}\left(|v_{j,n}|>\frac{\sqrt{\varepsilon p_j}}{\sqrt{6}n^{\frac{1}{4}}}\right)
+\mathbf{P}\left(|\hat{q}_{j,n}-q_j|>\frac{\sqrt{\varepsilon q_j}}{\sqrt{6}n^{\frac{1}{4}}}\right)
+\mathbf{P}\left(|v_{j,n}|>\frac{\sqrt{\varepsilon q_j}}{\sqrt{6}n^{\frac{1}{4}}}\right)
+\mathbf{P}\left(\left|v_{j,n}\right|>\frac{1}{2}\right)
$$
\beq\label{27.11.7}
\leq4\exp\{-C\sqrt{n}\}+\exp\{-Cn\}.
\eeq

Using Lemma~\ref{l.1} (see (\ref{25.11.3}), (\ref{25.11.5})) and Lemma~\ref{l.2}, for any $\varepsilon>0$, $1\leq j \leq r$ and some constant
$C:=C(\varepsilon,p_1,\dots,p_r,q_1,\dots,q_r,p,q)>0$ we obtain
$$
\mathbf{P}\left(\sqrt{n}\left|((\hat{p}_{j,n}-p_j)-(\hat{q}_{j,n}-q_j))\ln\frac{\hat{p}_{j,n}}{\hat{q}_{j,n}}-((\hat{p}_{j,n}-p_j)-(\hat{q}_{j,n}-q_j))\ln\frac{p_{j}}{q_{j}}\right|
\geq\varepsilon\right)
$$
$$
\leq\mathbf{P}\left(|\hat{p}_{j,n}-p_j|\cdot\left|\ln\frac{\hat{p}_{j,n}q_j}{p_j\hat{q}_{j,n}}\right|
\geq\frac{\varepsilon}{2\sqrt{n}}\right)
+\mathbf{P}\left(|\hat{q}_{j,n}-q_j|\cdot\left|\ln\frac{\hat{p}_{j,n}q_j}{p_j\hat{q}_{j,n}}\right|
\geq\frac{\varepsilon}{2\sqrt{n}}\right)
$$
$$
\leq\mathbf{P}\left(|\hat{p}_{j,n}-p_j|
\geq\frac{\sqrt{\varepsilon}}{\sqrt{2}n^{\frac{1}{4}}}\right)
+2\mathbf{P}\left(\left|\ln\frac{\hat{p}_{j,n}q_j}{p_j\hat{q}_{j,n}}\right|
\geq\frac{\sqrt{\varepsilon}}{\sqrt{2}n^{\frac{1}{4}}}\right)
$$
\beq\label{27.11.8}
+\mathbf{P}\left(|\hat{q}_{j,n}-q_j|
\geq\frac{\sqrt{\varepsilon}}{\sqrt{2}n^{\frac{1}{4}}}\right)
\leq 4\exp\{-C\sqrt{n}\}.
\eeq
Denote
$$
\xi_n:=\sqrt{n}\sum\limits_{j=1}^r\left((\hat{p}_{j,n}-p_j)\left(1+\ln\frac{p_j}{q_j}-\frac{q_j}{p_j}\right)
+(\hat{q}_{j,n}-q_j))\left(1+\ln\frac{q_j}{p_j}-\frac{p_j}{q_j}\right)\right).
$$
Relations (\ref{27.11.11}), (\ref{26.11.1}), (\ref{27.11.7}), (\ref{27.11.8}) imply that
for any $\varepsilon>0$ there exist constants $C_1:=C_1(\varepsilon,p_1,\dots,p_r,q_1,\dots,q_r,p,q,r)>0$,
$C:=C(\varepsilon,p_1,\dots,p_r,q_1,\dots,q_r,p,q)>0$
such that 
$$
\mathbf{P}(|\eta_n-\xi_n|>\varepsilon)\leq C_1\exp\{-C\sqrt{n}\}.
$$
Thus, thanks the Slutsky's theorem the random sequences $\eta_n$ and $\xi_n$ have the same limit by distribution. Finally, applying Lemma~\ref{l.3} with
$$
b_j=1+\ln\frac{p_j}{q_j}-\frac{q_j}{p_j}, \ \ \ c_j=1+\ln\frac{q_j}{p_j}-\frac{p_j}{q_j}
$$
we finish the proof.$\Box$

\end{section}

\begin{section}{Auxiliary results}

\begin{Lemma} \label{l.1} For any $g>0$ the following inequalities hold
\beq\label{25.11.1}
\mathbf{P}\left(\left|\hat{p}_n-p)\right|>g\right)\leq 2 \exp\left\{-\frac{n g^2}{2(\max(p,q))^2}\right\},
\eeq
\beq\label{25.11.2}
\mathbf{P}\left(\left|\hat{q}_n-q)\right|>g\right)\leq 2 \exp\left\{-\frac{n g^2}{2(\max(p,q))^2}\right\},
\eeq
\beq\label{26.11.5}
\max\limits_{1\leq j \leq r}\mathbf{P}\left(\frac{1}{n}\sum\limits_{i=1}^n(\mathbf{I}(X_i=a_j,Y_i=1)-pp_j)>g\right)
\leq 2 \exp\left\{-\frac{n g^2}{2(\max(pp_{\max},1-pp_{\min}))^2}\right\},
\eeq
\beq\label{26.11.7}
\max\limits_{1\leq j \leq r}\mathbf{P}\left(\frac{1}{n}\sum\limits_{i=1}^n(\mathbf{I}(X_i=a_j,Y_i=1)-qq_j)>g\right)
\leq 2 \exp\left\{-\frac{n g^2}{2(\max(qq_{\max},1-qq_{\min}))^2}\right\},
\eeq
$$
\max\limits_{1\leq j \leq r}\mathbf{P}\left(\left|\hat{p}_{j,n}-p_j\right|>g\right)\leq2 \exp\left\{-\frac{n g^2p^2}{2^7(\max(p,q))^2p_{\max}^2}\right\}+2\exp\left\{-\frac{n g^2p^2}{8(\max(pp_{\max},1-pp_{\min}))^2}\right\}
$$
\beq\label{25.11.3}
+\exp\left\{-\frac{n p^2p_{\min}^2}{2(\max(pp_{\max},1-pp_{\min}))^2}\right\}+\exp\left\{-\frac{n p^2}{8(\max(p,q))^2}\right\},
\eeq
$$
\max\limits_{1\leq j \leq r}\mathbf{P}\left(\left|\hat{q}_{j,n}-q_j\right|>g\right)\leq 2 \exp\left\{-\frac{n g^2q^2}{2^7(\max(p,q))^2q_{\max}^2}\right\}+2\exp\left\{-\frac{n g^2q^2}{8(\max(qq_{\max},1-qq_{\min}))^2}\right\}
$$
\beq\label{25.11.5}
+\exp\left\{-\frac{n q^2q_{\min}^2}{2(\max(qq_{\max},1-qq_{\min}))^2}\right\}+\exp\left\{-\frac{n q^2}{8(\max(p,q))^2}\right\},
\eeq
where $p_{\min}:=\min\limits_{1\leq j \leq r}p_j$, $p_{\max}:=\min\limits_{1\leq j \leq r}p_j$,
$q_{\min}:=\min\limits_{1\leq j \leq r}p_j$, $q_{\max}:=\min\limits_{1\leq j \leq r}p_j$.
\end{Lemma}
\DD. In order to prove (\ref{25.11.1}) we use the Hoeffding inequality. Thus we obtain
$$
\mathbf{P}\left(\left|\hat{p}_n-p)\right|>g\right)=
\mathbf{P}\left(\left|\frac{1}{n}\sum\limits_{i=1}^n(\mathbf{I}(Y_i=1)-p)\right|>g\right)
$$
$$
\leq\mathbf{P}\left(\sum\limits_{i=1}^n(\mathbf{I}(Y_i=1)-p)> n g\right)
$$
$$
+\mathbf{P}\left(-\sum\limits_{i=1}^n(\mathbf{I}(Y_i=1)-p)> n g\right)
\leq 2 \exp\left\{-\frac{n g^2}{2(\max(p,q))^2}\right\}.
$$
Inequalities (\ref{25.11.2})--(\ref{26.11.7}) we obtain by the same way, therefore we omit their proof

Let us prove (\ref{25.11.3}). We have
$$
\mathbf{P}\left(\left|\hat{p}_{j,n}-p_j\right|>g\right)
=\mathbf{P}\left(\left|\frac{\frac{1}{n}\sum\limits_{i=1}^n\mathbf{I}(X_i=a_j,Y_i=1)}{\frac{1}{n}\sum\limits_{i=1}^n\mathbf{I}(Y_i=1)}
\pm\frac{\frac{1}{n}\sum\limits_{i=1}^n\mathbf{I}(X_i=a_j,Y_i=1)}{p}-p_j\right|>g\right)
$$
$$
\leq\mathbf{P}\left(\left|\frac{\frac{1}{n}\sum\limits_{i=1}^n\mathbf{I}(X_i=a_j,Y_i=1)}{\frac{1}{n}\sum\limits_{i=1}^n\mathbf{I}(Y_i=1)}
-\frac{\frac{1}{n}\sum\limits_{i=1}^n\mathbf{I}(X_i=a_j,Y_i=1)}{p}\right|>\frac{g}{2}\right)
$$
\beq\label{25.11.14}
+\mathbf{P}\left(\left|\frac{\frac{1}{n}\sum\limits_{i=1}^n\mathbf{I}(X_i=a_j,Y_i=1)}{p}-p_j\right|>\frac{g}{2}\right)
=:\mathbf{P}_1+\mathbf{P}_2.
\eeq
Denote
$$
A:=\left\{\frac{1}{n}\sum\limits_{i=1}^n\mathbf{I}(X_i=a_j,Y_i=1)>2pp_j\right\}, \ \ \
B:=\left\{\frac{1}{n}\sum\limits_{i=1}^n\mathbf{I}(Y_i=1)<\frac{p}{2}\right\}.
$$
Let us upper bound the $\mathbf{P}_1$. We have
$$
\mathbf{P}_1=
\mathbf{P}\left(\left|
\frac{\left(p-\frac{1}{n}\sum\limits_{i=1}^n\mathbf{I}(Y_i=1)\right)\frac{1}{n}\sum\limits_{i=1}^n\mathbf{I}(X_i=a_j,Y_i=1)}
{p\frac{1}{n}\sum\limits_{i=1}^n\mathbf{I}(Y_i=1)}\right|>\frac{g}{2}\right)
$$
$$
=\mathbf{P}\left(\left|p-\frac{1}{n}\sum\limits_{i=1}^n\mathbf{I}(Y_i=1)\right|
>\frac{gp\frac{1}{n}\sum\limits_{i=1}^n\mathbf{I}(Y_i=1)}{2\frac{1}{n}\sum\limits_{i=1}^n\mathbf{I}(X_i=a_j,Y_i=1)},\overline{A},\overline{B}\right)
+\mathbf{P}(A)+\mathbf{P}(B)
$$
\beq\label{25.11.12}
\leq\mathbf{P}\left(\left|\frac{1}{n}\sum\limits_{i=1}^n(\mathbf{I}(Y_i=1)-p)\right|
>\frac{gp}{8p_j}\right)+\mathbf{P}(A)+\mathbf{P}(B).
\eeq
Let us upper bound the first term in the righthand side. Utilizing (\ref{25.11.1}), for any $1\leq j \leq r$ we have
\beq\label{25.11.7}
\mathbf{P}\left(\left|\frac{1}{n}\sum\limits_{i=1}^n(\mathbf{I}(Y_i=1)-p)\right|
>\frac{gp}{8p_j}\right)\leq2 \exp\Biggl\{-\frac{n g^2p^2}{2^7(\max(p,q))^2\Big(\max\limits_{1\leq j \leq r}p_j\Big)^2}\Biggl\}.
\eeq
Hoeffding inequality provide the following upper bound for the second term
$$
\mathbf{P}(A)=\mathbf{P}\left(\sum\limits_{i=1}^n(\mathbf{I}(X_i=a_j,Y_i=1)-pp_j)>npp_j\right)
$$
\beq\label{25.11.8}
\leq \exp\left\{-\frac{n p^2\left(\min\limits_{1\leq j \leq r}p_j\right)^2}{2\left(\max\limits_{1\leq j \leq r}(\max(pp_j,1-pp_j))\right)^2}\right\}.
\eeq
Now, let us upper bound $\mathbf{P}(B)$.Using again the Hoeffding inequality we obtain
\beq\label{25.11.10}
\mathbf{P}(B)=\mathbf{P}\left(\sum\limits_{i=1}^n(p-\mathbf{I}(Y_i=1))>\frac{p n}{2}\right)
\leq\exp\left\{-\frac{n p^2}{8(\max(p,q))^2}\right\}.
\eeq
Applying Hoeffding inequality for $\mathbf{P}_2$ we obtain
\beq\label{25.11.11}
\mathbf{P}_2\leq 2\exp\left\{-\frac{n g^2p^2}{8\left(\max\limits_{1\leq j \leq r}(\max(pp_j,1-pp_j))\right)^2}\right\}.
\eeq
Inequalities (\ref{25.11.14})--(\ref{25.11.11}) imply the relation (\ref{25.11.3}). In the same manner we can prove the inequality (\ref{25.11.5}). $\Box$

\begin{Lemma} \label{l.2} For any $g>0$ the following inequality holds
$$
\max\limits_{1\leq j \leq r}\mathbf{P}\left(\left|\ln\left(\frac{\hat{p}_{j,n}q_j}{p_j\hat{q}_{j,n}}\right)\right|>g\right)\leq4 \exp\left\{-\frac{n g^2p^2p^2_{\min}}{2^{11}(\max(p,q))^2p_{\max}^2}\right\}
$$
$$
+4\exp\left\{-\frac{n g^2p^2p^2_{\min}}{2^7(\max(pp_{\max},1-pp_{\min}))^2}\right\}+
2 \exp\left\{-\frac{n p^2p^2_{\min}}{2^{9}(\max(p,q))^2p_{\max}^2}\right\}
$$
$$
+2\exp\left\{-\frac{n p^2p^2_{\min}}{2^5(\max(pp_{\max},1-pp_{\min}))^2}\right\}
+3\exp\left\{-\frac{n p^2p_{\min}^2}{2(\max(pp_{\max},1-pp_{\min}))^2}\right\}
$$
$$
+3\exp\left\{-\frac{n p^2}{8(\max(p,q))^2}\right\}
$$
$$
+4 \exp\left\{-\frac{n g^2q^2q^2_{\min}}{2^{11}(\max(p,q))^2q_{\max}^2}\right\}+4\exp\left\{-\frac{n g^2q^2q^2_{\min}}{2^7(\max(qq_{\max},1-qq_{\min}))^2}\right\}
$$
$$
+2 \exp\left\{-\frac{n q^2q^2_{\min}}{2^{9}(\max(p,q))^2q_{\max}^2}\right\}+2\exp\left\{-\frac{n q^2q^2_{\min}}{2^5(\max(qq_{\max},1-qq_{\min}))^2}\right\}
$$
\beq\label{25.11.18}
+3\exp\left\{-\frac{n q^2q_{\min}^2}{2(\max(qq_{\max},1-qq_{\min}))^2}\right\}
+3\exp\left\{-\frac{n q^2}{8(\max(p,q))^2}\right\}.
\eeq
\end{Lemma}
\DD. For any $g>0$ and $1\leq j \leq r$ we obtain
\beq\label{25.11.15}
\mathbf{P}\left(\left|\ln\left(\frac{\hat{p}_{j,n}q_j}{p_j\hat{q}_{j,n}}\right)\right|>g\right)
=\mathbf{P}\left(\ln\left(\frac{\hat{p}_{j,n}q_j}{p_j\hat{q}_{j,n}}\right)>g\right)+
\mathbf{P}\left(-\ln\left(\frac{\hat{p}_{j,n}q_j}{p_j\hat{q}_{j,n}}\right)>g\right)=:\mathbf{P}_1+\mathbf{P}_2.
\eeq
Let us bound $\mathbf{P}_1$ from above. Using inequality $e^g-1\geq g$ and using Lemma~\ref{l.1}, for any $1\leq j \leq r$ we obtain
$$
\mathbf{P}_1=
\mathbf{P}\left(\frac{\hat{p}_{j,n}q_j}{p_j\hat{q}_{j,n}}>e^{g}\right)=
\mathbf{P}\left(\frac{\hat{p}_{j,n}q_j}{p_j\hat{q}_{j,n}}-1>e^{g}-1\right)
\leq\mathbf{P}\left(\frac{\hat{p}_{j,n}q_j}{p_j\hat{q}_{j,n}}-1>g\right)
$$
$$
\leq\mathbf{P}\left(\frac{\hat{p}_{j,n}q_j-p_j\hat{q}_{j,n}}{p_j\hat{q}_{j,n}}>g,\hat{q}_{j,n}\geq\frac{q_j}{2}\right)
+\mathbf{P}\left(\hat{q}_{j,n}<\frac{q_j}{2}\right)\leq
\mathbf{P}\left(\hat{p}_{j,n}q_j-p_j\hat{q}_{j,n}>\frac{g_1p_jq_j}{2}\right)
$$
$$
+\mathbf{P}\left(\hat{q}_{j,n}<\frac{q_j}{2}\right)
=\mathbf{P}\left(q_j(\hat{p}_{j,n}-p_j)-p_j(\hat{q}_{j,n}-q_j)>\frac{gp_jq_j}{2}\right)
+\mathbf{P}\left(\hat{q}_{j,n}<\frac{q_j}{2}\right)
$$
$$
\leq\mathbf{P}\left(q_j|\hat{p}_{j,n}-p_j|+p_j|\hat{q}_{j,n}-q_j|>\frac{gp_jq_j}{2}\right)
+\mathbf{P}\left(\hat{q}_{j,n}<\frac{q_j}{2}\right)
$$
$$
\leq\mathbf{P}\left(|\hat{p}_{j,n}-p_j|>\frac{gp_j}{4}\right)
+\mathbf{P}\left(|\hat{q}_{j,n}-q_j|>\frac{gq_j}{4}\right)
+\mathbf{P}\left(\hat{q}_{j,n}<\frac{q_j}{2}\right)
$$
$$
=\mathbf{P}\left(|\hat{p}_{j,n}-p_j|>\frac{gp_j}{4}\right)
+\mathbf{P}\left(|\hat{q}_{j,n}-q_j|>\frac{gq_j}{4}\right)
+\mathbf{P}\left(q_j-\hat{q}_{j,n}>\frac{q_j}{2}\right)
$$
$$
\leq2 \exp\left\{-\frac{n g^2p^2p^2_{\min}}{2^{11}(\max(p,q))^2p_{\max}^2}\right\}+2\exp\left\{-\frac{n g^2p^2p^2_{\min}}{2^7(\max(pp_{\max},1-pp_{\min}))^2}\right\}
$$
$$
+\exp\left\{-\frac{n p^2p_{\min}^2}{2(\max(pp_{\max},1-pp_{\min}))^2}\right\}+\exp\left\{-\frac{n p^2}{8(\max(p,q))^2}\right\}
$$
$$
+2 \exp\left\{-\frac{n g^2q^2q^2_{\min}}{2^{11}(\max(p,q))^2q_{\max}^2}\right\}+2\exp\left\{-\frac{n g^2q^2q^2_{\min}}{2^7(\max(qq_{\max},1-qq_{\min}))^2}\right\}
$$
$$
+2 \exp\left\{-\frac{n q^2q^2_{\min}}{2^{9}(\max(p,q))^2q_{\max}^2}\right\}+2\exp\left\{-\frac{n q^2q^2_{\min}}{2^5(\max(qq_{\max},1-qq_{\min}))^2}\right\}
$$
\beq\label{25.11.16}
+2\exp\left\{-\frac{n q^2q_{\min}^2}{2(\max(qq_{\max},1-qq_{\min}))^2}\right\}
+2\exp\left\{-\frac{n q^2}{8(\max(p,q))^2}\right\}.
\eeq
In the same way we can obtain the upper bound for  $\mathbf{P}_2$
$$
\mathbf{P}_2
\leq2 \exp\left\{-\frac{n g^2p^2p^2_{\min}}{2^{11}(\max(p,q))^2p_{\max}^2}\right\}+2\exp\left\{-\frac{n g^2p^2p^2_{\min}}{2^7(\max(pp_{\max},1-pp_{\min}))^2}\right\}
$$
$$
2 \exp\left\{-\frac{n p^2p^2_{\min}}{2^{9}(\max(p,q))^2p_{\max}^2}\right\}+2\exp\left\{-\frac{n p^2p^2_{\min}}{2^5(\max(pp_{\max},1-pp_{\min}))^2}\right\}
$$
$$
+2\exp\left\{-\frac{n p^2p_{\min}^2}{2(\max(pp_{\max},1-pp_{\min}))^2}\right\}+2\exp\left\{-\frac{n p^2}{8(\max(p,q))^2}\right\}
$$
$$
+2 \exp\left\{-\frac{n g^2q^2q^2_{\min}}{2^{11}(\max(p,q))^2q_{\max}^2}\right\}+2\exp\left\{-\frac{n g^2q^2q^2_{\min}}{2^7(\max(qq_{\max},1-qq_{\min}))^2}\right\}
$$
\beq\label{25.11.17}
+\exp\left\{-\frac{n q^2q_{\min}^2}{2(\max(qq_{\max},1-qq_{\min}))^2}\right\}
+\exp\left\{-\frac{n q^2}{8(\max(p,q))^2}\right\}.
\eeq
Inequalities (\ref{25.11.15})--(\ref{25.11.17}) imply
(\ref{25.11.18}).$\Box$

\begin{Lemma} \label{l.3} The following convergence takes place
$$
\sqrt{n}\sum\limits_{j=1}^r\left((\hat{p}_{j,n}-p_j)b_j
+(\hat{q}_{j,n}-q_j)c_j\right)\stackrel{d}{\rightarrow}\xi\sim\Phi_{0,\sigma^2},
  \ \text{as} \ n\rightarrow\infty,
$$
where
$$
\sigma^2=\mathbf{E}\Biggl(
\sum\limits_{j=1}^r\Biggl(\left(\frac{1}{p}\sum\limits_{i=1}^n(\mathbf{I}(X_i=a_j,Y_i=1)-pp_j)
-p_j\sum\limits_{i=1}^n(\mathbf{I}(Y_i=1)-p)\right)b_j
$$
$$
+\left(\frac{1}{q}\sum\limits_{i=1}^n(\mathbf{I}(X_i=a_j,Y_i=0)-qq_j)
-q_j\sum\limits_{i=1}^n(\mathbf{I}(Y_i=0)-q)\right)c_j\Biggl)\Biggl)^2,
$$
$\stackrel{d}{\rightarrow}$ means the convergence by the distribution, $\Phi_{0,\sigma^2}$ means the normal distributionwith mean $0$ and variance $\sigma^2$.
\end{Lemma}
\DD. Denote
$$
y_{j,n}:=
\frac{1}{\sqrt{n}p}\sum\limits_{i=1}^n(\mathbf{I}(X_i=a_j,Y_i=1)-pp_j)-\frac{p_j}{\sqrt{n}}\sum\limits_{i=1}^n(\mathbf{I}(Y_i=1)-p),
$$
$$
z_n:=\frac{1}{np}\sum\limits_{i=1}^n(\mathbf{I}(Y_i=1)-p),
$$
$$
v_{j,n}:=\frac{1}{\sqrt{n}p}\sum\limits_{i=1}^n(\mathbf{I}(X_i=a_j,Y_i=1)-pp_j).
$$
We have
$$
|\sqrt{n}(\hat{p}_{j,n}-p_j)-y_{j,n}|
=\left|\frac{\frac{1}{\sqrt{n}}\sum\limits_{i=1}^n\mathbf{I}(X_i=a_j,Y_i=1)}
{p\left(1+\frac{1}{np}\sum\limits_{i=1}^n(\mathbf{I}(Y_i=1)-p)\right)}-\sqrt{n}p_j-y_{j,n}\right|
$$
$$
=\left|\frac{1}{\sqrt{n}p}\sum\limits_{i=1}^n\mathbf{I}(X_i=a_j,Y_i=1)
\left(1-z_n+\sum\limits_{k=2}^\infty(-z_n)^k\right)-\sqrt{n}p_j-y_{j,n}\right|\mathbf{I}\left(|z_n|\leq\frac{1}{2}\right)
$$
$$
+|\sqrt{n}(\hat{p}_{j,n}-p_j)-y_{j,n}|\mathbf{I}\left(|z_n|>\frac{1}{2}\right)
$$
$$
=\left|(v_{j,n}+\sqrt{n}p_j)
\left(1-z_n+\sum\limits_{k=2}^\infty(-z_n)^k\right)-\sqrt{n}p_j-y_{j,n}\right|\mathbf{I}\left(|z_n|\leq\frac{1}{2}\right)
$$
$$
+|\sqrt{n}(\hat{p}_{j,n}-p_j)-y_{j,n}|\mathbf{I}\left(|z_n|>\frac{1}{2}\right)
$$
$$
=\left|-v_{j,n}z_n+
\frac{1}{\sqrt{n}p}\sum\limits_{i=1}^n\mathbf{I}(X_i=a_j,Y_i=1)
\sum\limits_{k=2}^\infty(-z_n)^k\right|\mathbf{I}\left(|z_n|\leq\frac{1}{2}\right)
$$
$$
+|\sqrt{n}(\hat{p}_{j,n}-p_j)-y_{j,n}|\mathbf{I}\left(|z_n|>\frac{1}{2}\right)
$$
$$
\leq|v_{j,n}z_n|\mathbf{I}\left(|z_n|\leq\frac{1}{2}\right)+
\left|
\frac{2}{\sqrt{n}p}\sum\limits_{i=1}^n\mathbf{I}(X_i=a_j,Y_i=1)z_n^2\right|\mathbf{I}\left(|z_n|\leq\frac{1}{2}\right)
$$
\beq\label{26.11.8}
+|\sqrt{n}(\hat{p}_{j,n}-p_j)-y_{j,n}|\mathbf{I}\left(|z_n|>\frac{1}{2}\right) \ \ \ \text{a.s.}
\eeq
Utilizing (\ref{26.11.8}) for any $1\leq j \leq r$, $\varepsilon>0$ we obtain
$$
\mathbf{P}(|\sqrt{n}(\hat{p}_{j,n}-p_j)-y_{j,n}|>\varepsilon)
\leq\mathbf{P}\left(|v_{j,n}z_n|\mathbf{I}\left(|z_n|\leq\frac{1}{2}\right)>\frac{\varepsilon}{3}\right)
$$
$$
+\mathbf{P}\left(\left|
\frac{2}{\sqrt{n}p}\sum\limits_{i=1}^n\mathbf{I}(X_i=a_j,Y_i=1)z_n^2\right|\mathbf{I}\left(|z_n|\leq\frac{1}{2
}\right)>\frac{\varepsilon}{3}\right)
$$
$$
+\mathbf{P}\left(|\sqrt{n}(\hat{p}_{j,n}-p_j)-y_{j,n}|\mathbf{I}\left(|z_n|>\frac{1}{2}\right)>\frac{\varepsilon}{3}\right)
$$
$$
\leq\mathbf{P}\left(|v_{j,n}z_n|>\frac{\varepsilon}{3}\right)
+\mathbf{P}\left(\left|
\frac{2}{\sqrt{n}p}\sum\limits_{i=1}^n\mathbf{I}(X_i=a_j,Y_i=1)z_n^2\right|>\frac{\varepsilon}{3}\right)
+\mathbf{P}\left(|z_n|>\frac{1}{2}\right)
$$
\beq\label{26.11.11}
=:\mathbf{P}_1+\mathbf{P}_2+\mathbf{P}_3.
\eeq
Let us upper bound the $\mathbf{P}_1$. Lemma~\ref{l.1} (see (\ref{25.11.1}), (\ref{26.11.5})) implies that for any $\varepsilon>0$ there exists 
$C:=C(\varepsilon,p_1,\dots,p_r,q_1,\dots,q_r,p,q)>0$ such that
$$
\mathbf{P}_1=\mathbf{P}\left(\left|\frac{1}{n}\sum\limits_{i=1}^n(\mathbf{I}(X_i=a_j,Y_i=1)- pp_j)\right|\cdot
\left|\frac{1}{n}\sum\limits_{i=1}^n(\mathbf{I}(Y_i=1)-p)
\right|>\frac{p^2\varepsilon}{3\sqrt{n}}\right)
$$
$$
\leq\mathbf{P}\left(\left|\frac{1}{n}\sum\limits_{i=1}^n(\mathbf{I}(X_i=a_j,Y_i=1)- pp_j)\right|>\frac{p\sqrt{\varepsilon}}{\sqrt{3}n^{\frac{1}{4}}}\right)
+\mathbf{P}\left(
\left|\frac{1}{n}\sum\limits_{i=1}^n(\mathbf{I}(Y_i=1)-p)
\right|>\frac{p\sqrt{\varepsilon}}{\sqrt{3}n^{\frac{1}{4}}}\right)
$$
\beq\label{26.11.12}
\leq2\exp\left\{-C\sqrt{n}\right\}.
\eeq
Let us upper bound the probability $\mathbf{P}_2$. Lemma~\ref{l.1} (see (\ref{25.11.1}), (\ref{26.11.5})) imply that for any $\varepsilon>0$ there exists
$C:=C(\varepsilon,p_1,\dots,p_r,q_1,\dots,q_r,p,q)>0$ such that
$$
\mathbf{P}_2=\mathbf{P}\left(\left(\frac{1}{n}\sum\limits_{i=1}^n\mathbf{I}(X_i=a_j,Y_i=1)\right)
\left|\frac{1}{n}\sum\limits_{i=1}^n(\mathbf{I}(Y_i=1)-p)
\right|^2>\frac{p^3\varepsilon}{6\sqrt{n}}\right)
$$
$$
\leq\mathbf{P}\left(
\left|\frac{1}{n}\sum\limits_{i=1}^n(\mathbf{I}(Y_i=1)-p)
\right|^2>\frac{p^2\varepsilon}{3\sqrt{n}p_j},\frac{1}{n}\sum\limits_{i=1}^n\mathbf{I}(X_i=a_j,Y_i=1)\geq\frac{pp_j}{2}\right)
$$
$$
+\mathbf{P}\left(\frac{1}{n}\sum\limits_{i=1}^n\mathbf{I}(X_i=a_j,Y_i=1)<\frac{pp_j}{2}\right)
$$
$$
\leq\mathbf{P}\left(
\left|\frac{1}{n}\sum\limits_{i=1}^n(\mathbf{I}(Y_i=1)-p)
\right|>\frac{p\sqrt{\varepsilon}}{\sqrt{3 p_j} n^{\frac{1}{4}}}\right)
+\mathbf{P}\left(\frac{1}{n}\left|\sum\limits_{i=1}^n (pp_j-\mathbf{I}(X_i=a_j,Y_i=1))\right|>\frac{pp_j}{2}\right)
$$
\beq\label{26.11.15}
\leq\exp\left\{-C\sqrt{n}\right\}+\exp\left\{-Cn\right\}.
\eeq
Utilizing Lemma~\ref{l.1} (see the relation (\ref{25.11.1})), we obtain the following: for any $\varepsilon>0$ and some  $C:=C(\varepsilon,p_1,\dots,p_r,q_1,\dots,q_r,p,q)>0$ we have
\beq\label{26.11.16}
\mathbf{P}_3=\mathbf{P}\left(
\left|\frac{1}{n}\sum\limits_{i=1}^n(\mathbf{I}(Y_i=1)-p)
\right|>\frac{1}{2}\right)\leq\exp\left\{-Cn\right\}.
\eeq
Relations (\ref{26.11.11})--(\ref{26.11.16}) imply that for any $\varepsilon>0$ there exists
$C:=C(\varepsilon,p_1,\dots,p_r,q_1,\dots,q_r,p,q)>0$ such that
\beq\label{26.11.18}
\max\limits_{1\leq j \leq r}\mathbf{P}(|\sqrt{n}(\hat{p}_{j,n}-p_j)-y_{j,n}|>\varepsilon)\leq 4\exp\left\{-C\sqrt{n}\right\}.
\eeq
Denote
$$
y_{j,n}':=
\frac{1}{\sqrt{n}q}\sum\limits_{i=1}^n(\mathbf{I}(X_i=a_j,Y_i=0)-qq_j)-\frac{q_j}{\sqrt{n}}\sum\limits_{i=1}^n(\mathbf{I}(Y_i=0)-q).
$$
In completely similar way as above for any $\varepsilon>0$ and some $C:=C(\varepsilon,p_1,\dots,p_r,q_1,\dots,q_r,p,q)>0$
we have
\beq\label{26.11.19}
\max\limits_{1\leq j \leq r}\mathbf{P}(|\sqrt{n}(\hat{q}_{j,n}-q_j)-y_{j,n}'|>\varepsilon)\leq 4\exp\left\{-C\sqrt{n}\right\}.
\eeq
It is easy to see that
$$
\sqrt{n}\sum\limits_{j=1}^r\left((\hat{p}_{j,n}-p_j)b_j
+(\hat{q}_{j,n}-q_j)c_j\right)=\sum\limits_{j=1}^r\left(y_jb_j
+y'_jc_j\right)+\sum\limits_{j=1}^r((\sqrt{n}(\hat{p}_{j,n}-p_j)-y_j)b_j)
$$
$$
+\sum\limits_{j=1}^r((\sqrt{n}(\hat{q}_{j,n}-q_j)-y_j')c_j)=:\xi_n+\eta_n+\zeta_n.
$$
By (\ref{26.11.18}), (\ref{26.11.19}) for any $\varepsilon>0$ we obtain
$$
\lim\limits_{n\rightarrow\infty}\mathbf{P}(|\eta_n+\zeta_n|>\varepsilon)=0.
$$
Therefore, by Slutsky's theorem, the weak limit of the original sequence coincides with the limit sequences $\xi_n$. It remains to note that
$$
\begin{aligned}
\xi_n=\frac{1}{\sqrt{n}}\sum\limits_{i=1}^n
\sum\limits_{j=1}^r\Biggl( & \left(\frac{1}{p}\sum\limits_{i=1}^n(\mathbf{I}(X_i=a_j,Y_i=1)-pp_j)
-p_j\sum\limits_{i=1}^n(\mathbf{I}(Y_i=1)-p)\right)b_j 
\\  &+\left(\frac{1}{q}\sum\limits_{i=1}^n(\mathbf{I}(X_i=a_j,Y_i=1)-qq_j)
-q_j\sum\limits_{i=1}^n(\mathbf{I}(Y_i=1)-q)\right)c_j \Biggl)
\end{aligned}
$$
and to apply central limit theorem. $\Box$

\end{section}

\section*{Acknowledgements}

We thanks Anatoly Yambartsev thank  for fruitful discussions.

\end{document}